\documentclass[11pt]{article}
\usepackage{amssymb,amsmath,booktabs,hyperref,mathtools}
\usepackage{graphics,authblk}
\usepackage{xcolor}
\usepackage{longtable,nicefrac}
\usepackage{graphicx}
\usepackage{subcaption}
\usepackage{multirow}
\usepackage{tabularx,enumerate}
\usepackage{rotating}
\usepackage{dcolumn,accents}
\usepackage{pdflscape}
\usepackage{array}
\usepackage{lscape}
\usepackage{anysize}
\usepackage{color}
\usepackage{amsthm}
\usepackage{listings}
\newcommand{\mathleft}{\@fleqntrue\@mathmargin0pt}
\newcommand{\mathcenter}{\@fleqnfalse}

\providecommand{\keywords}[1]{\textbf{\textit{Keywords: }} #1}
\usepackage{float}
\usepackage{soulutf8}
\usepackage{adjustbox}
\usepackage{colortbl}
\usepackage{blindtext,nameref}

\newtheorem{theorem}{Theorem}

\newtheorem{proposition}{Proposition}

\theoremstyle{definition}
\newtheorem{definition}{Definition}
\theoremstyle{remark}

\graphicspath{{photos/}}
\allowdisplaybreaks
\begin{document}
	\title{An analytic framework for the multiplicative best-worst method}
	\author{Harshit M. Ratandhara, Mohit Kumar}
	\date{}
	\affil{Department of Basic Sciences,\\ Institute of Infrastructure, Technology, Research and Management, Ahmedabad, Gujarat-380026, India\\ Email: harshitratandhara1999@gmail.com, msharmadma.iitr@gmail.com}
	\maketitle
	\begin{abstract}
		The Best-Worst Method (BWM) is a well-known Multi-Criteria Decision-Making (MCDM) method. This article deals with the multiplicative model of BWM. We first formulate an optimization model that is equivalent to the existing multiplicative model. This model provides a solid foundation for obtaining an analytic form of optimal interval-weights, Consistency Index (CI) and Consistency Ratio (CR). The proposed approach does not require any optimization software, which makes it easy to implement as well as time efficient. Also, the obtained analytical form of CR permits it to serve as an input-based consistency measure. After obtaining these analytic forms, a secondary objective function is introduced to select the best optimal weight set from the collection of all optimal weight sets. Finally, we discuss some numerical examples and a real-world application of the proposed approach in ranking the drivers of Sustainable Additive Manufacturing (SAM) to illustrate the proposed approach.
	\end{abstract}
	\section{Introduction}
	The discrete Multi-Criteria Decision-Making (MCDM), also known as Multi-Attribute Decision-Making (MADM), is one of the key branches of the operation research which mainly deals with the selection and the ranking problems. We shall refer to discrete MCDM simply as MCDM because these terms are often used interchangeably in the existing literature. The process of tackling an MCDM problem usually involves two steps: first, the computation of criteria-weights, and second, the ranking of alternatives. Based on these steps, MCDM methods can be broadly categorized into two classes, weight calculation methods and ranking methods. AHP\cite{saaty1990make}, SMARTS\cite{edwards1994smarts}, etc. are weight calculation methods and TOPSIS\cite{hwang1981methods}, VIKOR\cite{opricovic2004compromise}, ELECTRE\cite{roy1990outranking}, PROMETHEE\cite{brans1985note}, etc. are ranking methods.\\\\
	The pairwise comparison methods form an important subclass of weight calculation methods. These methods use a matrix called pairwise comparison matrix in calculating weights. Let $C=\{c_1,c_2,...,c_n\}$ be a set of decision-criteria. Then the pairwise comparison matrix is of the form $A=[a_{ij}]_{n\times n}$, where $a_{ij}$ represents the relative preference of $i^{th}$ criterion over $j^{th}$ criterion\cite{brunelli2018survey}. The pairwise comparison methods are divided into two groups: eigenvector methods and extremal methods. Eigenvector methods utilize the principal eigenvector of $A$, whereas extremal methods are distance based methods\cite{golany1993multicriteria}.\\\\
	The Best-Worst Method (BWM) is an extremal pairwise comparison method that uses maximum deviation as a distance function\cite{rezaei2015best}. The first step of BWM is the selection of two reference criteria, the best and the worst. Then the best-to-other vector $A_{bo}=(a_{b1},a_{b2},...,a_{bn})$ and the other-to-worst vector $A_{ow}=(a_{1w},a_{2w},...,a_{nw})$ are determined, which forms the $b^{th}$ row and the $w^{th}$ column of $A$, respectively. Using these comparison values, a non-linear minimization problem is formulated. An optimal solution of this problem gives an optimal weight set. Apart from this non-linear model, some other models of BWM using different distance functions such as Euclidean distance\cite{kocak2018euclidean}, taxicab distance (goal programming model)\cite{amiri2020goal}, etc. are also developed.
	Estimating the accuracy of a weight set is an essential part of any BWM-model. For that purpose, Rezaei\cite{rezaei2015best} proposed the concepts of Consistency Index (CI) and Consistency Ratio (CR) for the non-linear model, which are further extended for most of the other models as well. Some other significant theoretical developments in BWM are as follow. To deal with the non-uniqueness of optimal weight sets in the non-linear model, Rezaei\cite{rezaei2016best} calculated optimal interval-weights, and then employed interval-analysis to weight the criteria. He also proposed a linear model to obtain a unique weight set\cite{rezaei2016best}. Safarzadeh et al.\cite{safarzadeh2018group} introduced BWM in group decision-making. To analyze the optimal interval-weights of the non-linear model, Rezaei\cite{rezaei2020concentration} proposed a ratio called concentration ratio. Liang et al.\cite{liang2020consistency} introduced the concepts of input-based consistency ratio and ordinal consistency ratio for BWM. He also calculated threshold values to determine the degree upto which inconsistency in comparison values is acceptable\cite{liang2020consistency}. Mohammadi and Rezaei\cite{mohammadi2020bayesian} incorporated Bayesian with BWM to calculate the final aggregated weights. Liang et al.\cite{liang2021belief} proposed belief function based BWM to address uncertainty in preferences. Lei et al.\cite{lei2022preference} developed preference rationality analysis to eliminate ordinal inconsistencies in comparison values. Tu et al.\cite{tu2023priority} introduced two prioritization models of BWM along with their threshold values. Wu et al.\cite{wu2023analytical} proposed an analytic framework for the non-linear model. The BWM is also been extended for fuzzy set\cite{guo2017fuzzy, mohtashami2021novel, dong2021fuzzy}, intuitionistic fuzzy set\cite{wan2021novel, mou2016intuitionistic}, hesitant fuzzy set\cite{ali2019hesitant}, etc. The BWM has tremendous applications in various fields including supplier selection \cite{vahidi2018sustainable,ahmadi2017assessing}, automotive\cite{van2017battle}, eco-industrial parks\cite{zhao2018comprehensive}, airline industry\cite{rezaei2017multi,gupta2018evaluating}, energy efficiency\cite{gupta2017developing,wang2019energy}, location selection\cite{kheybari2020sustainable}, and many more. The BWM is also combined with some other MCDM methods such as BWM-TOPSIS\cite{youssef2020integrated}, BWM-VIKOR\cite{dawood2023novel}, BWM-MULTIMOORA\cite{sarabi2021developing}, BWM-SERVQUAL\cite{rezaei2018quality} and  Best-Worst Tradeoff method\cite{liang2022best}.\\\\
	Using the distance function $\max\{\nicefrac{x}{y},\nicefrac{y}{x}\}$ on $(\mathbb{R_+},\cdot,\leq)$, Brunelli and Rezaei\cite{brunelli2019multiplicative} developed a non-linear multiplicative model of BWM, which they linearized using the logarithmic transformation. It is known that this model may lead to multiple optimal weight sets. Although multiple weight sets are preferable in some cases, in most cases, the decision maker prefers a unique weight set. In this article, we mainly focus on two objectives: how to solve the multiplicative model analytically and how to obtain a unique weight set. To derive the analytical framework for the multiplicative model, we formulate a novel optimization model that transforms inconsistent pairwise comparison systems into consistent ones while preserving as much of the original information as possible. Such consistent systems are referred as optimally modified pairwise comparison systems. After establishing one-to-one correspondence between optimal weight sets and optimally modified pairwise comparison systems, we derive an analytic form of such systems, which ultimately leads to an analytic form of the optimal interval-weights. We also derive analytical forms of consistency index (CI) and consistency ratio (CR). This analytical framework is important in several aspects. It gives the optimal interval-weights directly in terms of the comparison values and thus, no optimization software is required in the implementation of the multiplicative model. The analytic form of CR, which is an input-based formulation, can be used to provide an immediate feedback to the decision maker regarding the consistency of the comparison values. To obtain the best optimal weight set, we follow the same path as Wu et al.\cite{wu2023analytical}. We introduce a secondary objective function to select the best optimally modified pairwise comparison system, which leads to the best optimal weight set. It is important to mention that secondary objective function retains all the features of the multiplicative model. Finally, we illustrate the proposed approach using numerical examples and a real-world application in ranking the drivers of sustainable additive manufacturing.\\\\
	The rest of the paper is organized as follows. Section 2 provides basic definitions and a brief introduction to multiplicative BWM. Section 3 presents an analytic framework for multiplicative BWM. A real-world application of the proposed model is discussed in Section 4, while Section 5 provides conclusive remarks and future directions. 
	\section{Preliminary}
	\subsection{Definition and Notations}
	\begin{definition}
		The best-to-other vector $A_{bo}$ and the other-to-worst vector $A_{ow}$ are jointly called Pairwise Comparison System (PCS), denoted as $(A_{bo},A_{ow})$.
	\end{definition}
	\begin{definition}\cite{rezaei2015best}
		A PCS $(A_{bo},A_{ow})$ is said to be consistent if $a_{bi}\times a_{iw}=a_{bw}$ for all $i\in\{1,2,...,n\}\setminus\{b,w\}$.
	\end{definition}
	Throughout the article, $C=\{c_1,c_2,...,c_n\}$ is the set of all decision-criteria and $D=\{c_1,c_2,...c_n\}\setminus\{c_b,c_w\}$. Whenever there is no ambiguity, we shall simply use $C=\{1,2,...,n\}$ and $D=\{1,2,...,n\}\setminus\{b,w\}$.\\\\
	If $(A_{bo},A_{ow})$ is consistent, \cite{wu2023analytical}
	\begin{equation}\label{weights}
		w_j=\frac{a_{jw}}{\sum_{i=1}^{n}a_{iw}}= \frac{1}{a_{bj}\sum_{i=1}^{n}\frac{1}{a_{bi}}}, \text{ for }j\in C,
	\end{equation} 
	is the unique solution of the system of equations
	\begin{equation}\label{system}
		\frac{w_b}{w_i}=a_{bi},\quad \frac{w_i}{w_w}=a_{iw}, \quad \frac{w_b}{w_w}=a_{bw}, \ i\in D.
	\end{equation}
	\subsection{Multiplicative Best-Worst Method}
	Brunelli and Rezaei\cite{brunelli2019multiplicative} developed BWM in the multiplicative framework, thus known as multiplicative BWM, using the $\mathcal{G}$-distance $d_\mathcal{G}$ for the Abelian linearly ordered group $(\mathbb{R_+},\cdot,\leq)$ defined as
	\begin{equation}\label{distance}
		d_\mathcal{G}(x,y)=\max\{\nicefrac{x}{y},\nicefrac{y}{x}\},\ x,y\in \mathbb{R}_+,
	\end{equation}
	where $(\mathbb{R}_+,\cdot,\leq)$ is the set of positive real numbers equipped with the usual multiplication and the usual order.\\\\
	In this model, after selecting the best and the worst criterion from $C$, the decision maker provides the relative preferences of the best criterion over each criterion and the relative preferences of each criterion over the worst criterion. These preferences are usually provides in the form of linguistic terms, which are assigned numerical values using some known scale such as $\frac{1}{9}$ to $9$ scale\cite{saaty1994make}. These preference values form $(A_{bo},A_{ow})$. Using this system, the following minimax optimization problem is formulated\cite{brunelli2019multiplicative}. 
	\begin{equation}\label{optimization_1}
		\begin{split}
			&\min\biggl\{ \max\{\nicefrac{a_{ij}}{\frac{w_i}{w_j}},\nicefrac{\frac{w_i}{w_j}}{a_{ij}}:(i,j)\in E \}\biggr\}\\
			&\text{subject to: }  w_1+w_2+...+w_n=1 \text{ and } w_k\geq0 \text{ for all }k=1,2,...,n,
		\end{split}
	\end{equation}	
	where $E=\{(i,j):(i=b \lor i=w)\land i\neq j\}$. An optimal solution of problem \eqref{optimization_1} is an optimal weight set.\\\\
	It is known that problem \eqref{optimization_1} may have multiple optimal solutions and hence we get multiple optimal weight sets\cite{brunelli2019multiplicative}. To deal with this problem, Brunelli and Rezaei\cite{brunelli2019multiplicative} linearized problem \eqref{optimization_1} using the transformation $b_{ij}=ln (a_{ij})$. In this linearization, the equation $v_1+v_2+...+v_n=0$, where $v_i=ln(w_i)$, is used as the normalization condition. This gives $ln(w_1)+ln(w_2)+...+ln(w_n)=0$, i.e., $w_1\times w_2\times...\times w_n=1$. Now, the fact that for any weight set $\{w_1,w_2,...,w_n\}$, there exists unique weight sets $\{w_1',w_2',...,w_n'\}$ and $\{w_1'',w_2'',...,w_n''\}$ such that $\frac{w_i}{w_j}=\frac{w_i'}{w_j'}=\frac{w_i''}{w_j''}$ for all $i,j$, $w_1'+w_2'+...+w_n'=1$ and $w_1''\times w_2''\times...\times w_n''=1$ implies that there is one-one correspondence between optimal solutions of problem \eqref{optimal_weights} and it's linearized model, but their feasible regions are different due to the difference in their normalization conditions. So, linearized model also gives multiple optimal weight sets. 
	\section{An analytic framework for the non-linear multiplicative best-worst method}
	In this section, we have three objectives: to obtain an analytic form of optimal weights, to introduce a secondary objective function to obtain a unique optimal weight set and to derive an analytic form of consistency index and consistency ratio.
	\subsection{An analytic form of optimal weights}
	Since the non-linear multiplicative model and linearized model are equivalent, it is sufficient to obtain an analytic form of optimal weights of the non-linear multiplicative model.\\\\
	Consider the following minimization problem.
	\begin{equation}\label{optimization_2}
		\begin{split}
			&\min\epsilon\\
			&\text{subject to:}\\
			&\nicefrac{a_{bi}}{\frac{w_b}{w_i}}\leq \epsilon, \quad \nicefrac{\frac{w_b}{w_i}}{a_{bi}}\leq \epsilon, \quad \nicefrac{a_{iw}}{\frac{w_i}{w_w}}\leq \epsilon, \quad \nicefrac{\frac{w_i}{w_w}}{a_{iw}}\leq \epsilon, \quad \nicefrac{a_{bw}}{\frac{w_b}{w_w}}\leq \epsilon, \quad \nicefrac{\frac{w_b}{w_w}}{a_{bw}}\leq \epsilon,\\
			&w_1+w_2+...+w_n=1, \quad w_k\geq0 \text{ for all } i\in D \text{ and } k\in C.
		\end{split}
	\end{equation}
	Note that problem \eqref{optimization_2} is an equivalent formulation of problem \eqref{optimization_1}. It is a non-linear problem having $n+1$ variables $w_1,w_2,...,w_n$ and $\epsilon$. So, it has optimal solution(s) of the form $(w_1^*,w_2^*,...,w_n^*,\epsilon^*)$. For each optimal solution, $W^*=\{w_1^*,w_2^*,...,w_n^*\}$ is an optimal weight set and $\epsilon^*$ represents the accuracy of $W^*$. Since $\epsilon^*$ is also the optimal objective value, it remains same for all $W^*$.\\\\
	Now consider another minimization problem.
	\begin{equation}\label{optimization_6}
		\begin{split}
			&\min \eta\\
			&\text{sub to:} \\
			&\nicefrac{a_{bi}}{\tilde{a}_{bi}}\leq \eta,\quad \nicefrac{\tilde{a}_{bi}}{a_{bi}}\leq \eta,\quad \nicefrac{a_{iw}}{\tilde{a}_{iw}}\leq \eta,\quad \nicefrac{\tilde{a}_{iw}}{a_{iw}}\leq \eta,\quad \nicefrac{a_{bw}}{\tilde{a}_{bw}}\leq \eta,\quad \nicefrac{\tilde{a}_{bw}}{a_{bw}}\leq\eta,\\
			&\tilde{a}_{bi}\times \tilde{a}_{iw}=\tilde{a}_{bw},\quad \tilde{a}_{bi},\tilde{a}_{iw},\tilde{a}_{bw}\geq 0 \text{ for all } i\in D.
		\end{split}
	\end{equation}
	Observe that problem \eqref{optimization_6} is a non-linear problem with $2n-2$ variables $\tilde{a}_{b1},...,\tilde{a}_{bb-1},\tilde{a}_{bb+1},...,\tilde{a}_{bn},$ $\tilde{a}_{1w},...,\tilde{a}_{w-1w},\tilde{a}_{w+1w},...,\tilde{a}_{nw}$ and $\eta$. So, it has optimal solution(s) of the form \\$(\tilde{a}_{b1}^*,...,\tilde{a}_{bb-1}^*,\tilde{a}_{bb+1}^*,...,\tilde{a}_{bn}^*,\tilde{a}_{1w}^*,...,\tilde{a}_{w-1w}^*,\tilde{a}_{w+1w}^*,...,\tilde{a}_{nw}^*,\eta^*)$. Each optimal solution along with $\tilde{a}_{bb}^*=\tilde{a}_{ww}^*=1$ forms a consistent PCS $(\tilde{A}_{bo}^*,\tilde{A}_{ow}^*)$, which is known as optimally modified PCS and $\eta^*$ is the measurement of accuracy of this PCS. Since $\eta^*$ is also the optimal objective value, it remains same for all $(\tilde{A}_{bo}^*,\tilde{A}_{ow}^*)$.\\\\
	Now, we establish the relationship between optimal solutions of problem \eqref{optimization_2} and problem \eqref{optimization_6}.\\\\
	Let $\{w_1^*,w_2^*,...,w_n^*\}$ be an optimal weight set. So, we have 
	\begin{equation}\label{equivalence} 
		\nicefrac{a_{bi}}{\frac{w_b^*}{w_i^*}}\leq \epsilon^*, \quad \nicefrac{\frac{w_b^*}{w_i^*}}{a_{bi}}\leq \epsilon^*, \quad \nicefrac{a_{iw}}{\frac{w_i^*}{w_w^*}}\leq \epsilon^*, \quad \nicefrac{\frac{w_i^*}{w_w^*}}{a_{iw}}\leq \epsilon^*, \quad \nicefrac{a_{bw}}{\frac{w_b^*}{w_w^*}}\leq \epsilon^*, \quad \nicefrac{\frac{w_b^*}{w_w^*}}{a_{bw}}\leq \epsilon^*
	\end{equation}
	for all $i\in D$.\\\\
	Take 
	\begin{equation}\label{optimal_pcs}
		\tilde{a}_{bi}=\frac{w_b^*}{w_i^*},\quad \tilde{a}_{iw}=\frac{w_i^*}{w_w^*}, \quad \tilde{a}_{bw}=\frac{w_b^*}{w_w^*}, \quad \tilde{a}_{bb}=\tilde{a}_{ww}=1, \quad i\in D.
	\end{equation}
	Observe that $(\tilde{A}_{bo},\tilde{A}_{ow})$ is a consistent PCS. Also,
	\begin{equation}\label{equivalence_1}
		\nicefrac{a_{bi}}{\tilde{a}_{bi}}\leq \epsilon^*,\quad \nicefrac{\tilde{a}_{bi}}{a_{bi}}\leq \epsilon^*,\quad \nicefrac{a_{iw}}{\tilde{a}_{iw}}\leq \epsilon^*,\quad \nicefrac{\tilde{a}_{iw}}{a_{iw}}\leq \epsilon^*,\quad \nicefrac{a_{bw}}{\tilde{a}_{bw}}\leq \epsilon^*,\quad \nicefrac{\tilde{a}_{bw}}{a_{bw}}\leq\epsilon^*
	\end{equation}
	for all $i\in D$. This gives $\eta^*\leq \epsilon^*$.\\\\
	Let $(\tilde{A}_{bo}^*,\tilde{A}_{ow}^*)$ be an optimally modified PCS. So, we have 
	\begin{equation}\label{equivalence_2}
		\nicefrac{a_{bi}}{\tilde{a}_{bi}^*}\leq \eta^*,\quad \nicefrac{\tilde{a}_{bi}^*}{a_{bi}}\leq \eta^*,\quad \nicefrac{a_{iw}}{\tilde{a}_{iw}^*}\leq \eta^*,\quad \nicefrac{\tilde{a}_{iw}^*}{a_{iw}}\leq \eta^*,\quad \nicefrac{a_{bw}}{\tilde{a}_{bw}^*}\leq \eta^*,\quad \nicefrac{\tilde{a}_{bw}^*}{a_{bw}}\leq\eta^*
	\end{equation}
	for all $i\in D$.\\\\
	By \eqref{weights}, the system of equations $\frac{w_b}{w_i}=\tilde{a}_{bi}^* $, $\frac{w_i}{w_w}=\tilde{a}_{iw}^*$ and $\frac{w_b}{w_w}=\tilde{a}_{bw}^*$, $i\in D$ has a unique solution
	\begin{equation}\label{optimal_weight}
		w_j=\frac{\tilde{a}_{jw}^*}{\sum_{i=1}^{n}\tilde{a}_{iw}^*},\quad j\in C.
	\end{equation} 
	This implies
	\begin{equation}\label{equivalence_3} 
		\nicefrac{a_{bi}}{\frac{w_b}{w_i}}\leq \eta^*, \quad \nicefrac{\frac{w_b}{w_i}}{a_{bi}}\leq \eta^*, \quad \nicefrac{a_{iw}}{\frac{w_i}{w_w}}\leq \eta^*, \quad \nicefrac{\frac{w_i}{w_w}}{a_{iw}}\leq \eta^*, \quad \nicefrac{a_{bw}}{\frac{w_b}{w_w}}\leq \eta^*, \quad \nicefrac{\frac{w_b}{w_w}}{a_{bw}}\leq \eta^*
	\end{equation}
	for all $i\in D$. This gives $\epsilon^*\leq \eta^*$.\\\\
	From the above discussion, we get $\epsilon^*=\eta^*$, i.e., problem \eqref{optimization_2} and problem \eqref{optimization_6} have the same optimal objective value. This implies that equation \eqref{optimal_pcs} and \eqref{optimal_weight} represents an optimally modified PCS and an optimal weight set respectively. So, for each $W^*=\{w_1^*,w_2^*,...,w_n^*\}$, there exists a unique $(\tilde{A}_{bo}^*,\tilde{A}_{ow}^*)$ such that $w_j^*=\frac{\tilde{a}_{jw}^*}{\sum_{i=1}^{n}\tilde{a}_{iw}^*}$. So, to get analytic solution(s) of problem \eqref{optimization_2}, it is sufficient to obtain analytic solution(s) of problem \eqref{optimization_6}. In this respect, we first introduce some notations.\\\\
	Let $(A_{bo},A_{ow})$ be a PCS. Then for $i,j\in D$, define
	\begin{equation}\label{CV}
		\epsilon_i= \max\biggl\{\frac{a_{bi}\times a_{iw}}{a_{bw}},\frac{a_{bw}}{a_{bi}\times a_{iw}}\biggr\}^{\frac{1}{3}} \quad\text{and} \quad 
		\epsilon_{i,j}= \max\biggl\{\frac{a_{bi}\times a_{iw}}{a_{bj}\times a_{jw}},\frac{a_{bj}\times a_{jw}}{a_{bi}\times a_{iw}}\biggr\}^{\frac{1}{4}}.
	\end{equation}
	Observe that equation \eqref{CV} can be rewritten as
	\begin{equation}\label{CV1}
		\epsilon_i= \begin{cases}
			\bigg(\frac{a_{bi}\times a_{iw}}{a_{bw}}\bigg)^{\frac{1}{3}}, \quad \text{if } a_{bi}\times a_{iw}\geq a_{bw}\\
			\bigg(\frac{a_{bw}}{a_{bi}\times a_{iw}}\bigg)^{\frac{1}{3}}, \quad \text{if } a_{bi}\times a_{iw}<a_{bw}
		\end{cases},
	\end{equation}
	\begin{equation}\label{CV2}
		\epsilon_{i,j}= \begin{cases}
			\bigg(\frac{a_{bi}\times a_{iw}}{a_{bj}\times a_{jw}}\bigg)^{\frac{1}{4}}, \quad \text{if } a_{bi}\times a_{iw}\geq a_{bj}\times a_{jw}\\
			\bigg(\frac{a_{bj}\times a_{jw}}{a_{bi}\times a_{iw}}\bigg)^{\frac{1}{4}}, \quad \text{if } a_{bi}\times a_{iw}<a_{bj}\times a_{jw}
		\end{cases}.
	\end{equation}
	Here, we consider only positive fourth root. Note that $\epsilon_i,\epsilon_{i,j}\geq 1$.
	\begin{proposition}\label{lower1}
		Let $\epsilon_j,\epsilon_{j,k}$ and $\eta^*$ be as above. Then $\epsilon_j,\epsilon_{j,k}\leq\eta^*$.
	\end{proposition}
	\begin{proof}
		Here, we prove only $\epsilon_j\leq \eta^*$ as the other part can be proven using the similar argument. Let $(\tilde{A}_{bo}^*,\tilde{A}_{ow}^*)$ be an optimally modified PCS. Let $\max\{\nicefrac{a_{bj}}{\tilde{a}_{bj}^*},\nicefrac{\tilde{a}_{bj}^*}{a_{bj}}\}=\eta_{bj}$, $\max\{\nicefrac{a_{jw}}{\tilde{a}_{jw}^*},\nicefrac{\tilde{a}_{jw}^*}{a_{jw}}\}=\eta_{jw}$ and $\max\{\nicefrac{a_{bw}}{\tilde{a}_{bw}^*},\nicefrac{\tilde{a}_{bw}^*}{a_{bw}}\}=\eta_{bw}$. Note that $\nicefrac{a_{bj}}{\tilde{a}_{bj}^*}\geq 1$ iff $\nicefrac{\tilde{a}_{bj}^*}{a_{bj}}\leq 1$ and vice versa. So, we get $1\leq\eta_{bj},\eta_{jw},\eta_{bw}\leq \eta^*$. Also, $\tilde{a}_{bj}^*\in\{\frac{a_{bj}}{\eta_{bj}},\eta_{bj}\times a_{bj}\}$, $\tilde{a}_{jw}^*\in\{\frac{a_{jw}}{\eta_{jw}},\eta_{jw}\times a_{jw}\}$ and $\tilde{a}_{bw}^*\in\{\frac{a_{bw}}{\eta_{bw}},\eta_{bw}\times a_{bw}\}$.\\\\
		If $a_{bj}\times a_{jw}=a_{bw}$, then by \eqref{CV1}, $\epsilon_j=1$ and we are done. If $a_{bj}\times a_{jw}>a_{bw}$, then by \eqref{CV1}, ${\epsilon_j}^3=\frac{a_{bj}\times a_{jw}}{a_{bw}}$. This gives $\frac{a_{bj}}{\epsilon_j}\times \frac{a_{jw}}{\epsilon_j}=\epsilon_j\times a_{bw}$. Now to prove $\epsilon_j\leq \eta^*$, it is sufficient to prove that at least one of $\epsilon_j\leq\eta_{bj}$, $\epsilon_j\leq\eta_{jw}$ and $\epsilon_j\leq\eta_{bw}$ holds. Suppose, if possible, neither of these conditions hold. Then we get $\frac{a_{bj}}{\epsilon_j}<\tilde{a}_{bj}^*$, $\frac{a_{jw}}{\epsilon_j}<\tilde{a}_{jw}^*$ and $\epsilon_j\times a_{bw}>\tilde{a}_{bw}^*$. This gives $\tilde{a}_{bw}^*<\tilde{a}_{bj}^*\times \tilde{a}_{jw}^*$, which is contradiction. Thus, $\epsilon_j\leq \eta^*$. Now, if $a_{bj}\times a_{jw}<a_{bw}$, then result can be proven by replicating the above argument. Hence the proof.
	\end{proof}
	Let $D_1=\{i\in D: a_{bi}\times a_{iw}<a_{bw}\}$ and $D_2=\{i\in D: a_{bi}\times a_{iw}>a_{bw}\}$. Let $i_0\in D_1$ and $j_0\in D_2$ be such that $a_{bi_0}\times a_{i_0w}\leq a_{bi}\times a_{iw}$ and $a_{bj_0}\times a_{j_0w}\geq a_{bi}\times a_{iw}$ for all $i\in D$. From Proposition \ref{lower1}, it follows that $\max \{\epsilon_{i_0},\epsilon_{j_0},\epsilon_{i_0,j_0}\}\leq \eta^*$. If $D_1 (D_2)=\phi$, then we simply ignore $i_0$ ($j_0$).
	\begin{theorem}\label{exact_obj}
		Let $i_0,j_0$ and $\eta^*$ be as above. Then the following statements hold.
		\begin{enumerate}
			\item If $\epsilon_{i_0}=\max \{\epsilon_{i_0},\epsilon_{j_0},\epsilon_{i_0,j_0}\}$, then $\eta^*=\epsilon_{i_0}$. Also, $(\tilde{A}_{bo},\tilde{A}_{ow})$ defined as
			\begin{equation}\label{m2}
				\tilde{a}_{bi}=\bigg(\frac{a_{bi}\times a_{bw}}{\epsilon_{i_0}\times a_{iw}}\bigg)^{\frac{1}{2}},\quad \tilde{a}_{iw}=\bigg(\frac{a_{iw}\times a_{bw}}{\epsilon_{i_0}\times a_{bi}}\bigg)^{\frac{1}{2}},\quad\tilde{a}_{bw}= \frac{a_{bw}}{\epsilon_{i_0}},\quad i\in D
			\end{equation}
			is an optimally modified PCS.
			\item If $\epsilon_{j_0}=\max \{\epsilon_{i_0},\epsilon_{j_0},\epsilon_{i_0,j_0}\}$, then $\eta^*=\epsilon_{j_0}$. Also, $(\tilde{A}_{bo},\tilde{A}_{ow})$ defined as
			\begin{equation}\label{m1}
				\tilde{a}_{bi}=\bigg(\frac{\epsilon_{j_0}\times a_{bi}\times a_{bw}}{a_{iw}}\bigg)^{\frac{1}{2}},\quad \tilde{a}_{iw}=\bigg(\frac{\epsilon_{j_0}\times a_{iw}\times  a_{bw}}{a_{bi}}\bigg)^{\frac{1}{2}},\quad\tilde{a}_{bw}= \epsilon_{j_0}\times a_{bw},\quad i\in D
			\end{equation}
			is an optimally modified PCS.
			\item If $\epsilon_{i_0,j_0                                                                                                                                                              }=\max \{\epsilon_{i_0},\epsilon_{j_0},\epsilon_{i_0,j_0}\}$, then $\eta^*=\epsilon_{i_0,j_0}$. Also, $(\tilde{A}_{bo},\tilde{A}_{ow})$ defined as
			\begin{equation}\label{m3}
				\begin{split}
					&\tilde{a}_{bi}=\bigg(\frac{\epsilon_{i_0,j_0}^2\times a_{bi}\times a_{bi_0}\times a_{i_0w}}{a_{iw}}\bigg)^{\frac{1}{2}},\quad \tilde{a}_{iw}=\bigg(\frac{\epsilon_{i_0,j_0}^2\times a_{iw}\times a_{bi_0}\times a_{i_0w}}{a_{bi}}\bigg)^{\frac{1}{2}},\\
					&\tilde{a}_{bw}= \epsilon_{i_0,j_0}^2\times a_{bi_0}\times a_{i_0w},\quad i\in D
				\end{split}
			\end{equation}
			is an optimally modified PCS. 
		\end{enumerate}
	\end{theorem}
	\begin{proof}
		Here, we shall prove only $1^{st}$ statement. Note that $(\tilde{A}_{bo},\tilde{A}_{ow})$ is consistent. Now, if we prove that
		\begin{equation*}
			\nicefrac{a_{bi}}{\tilde{a}_{bi}}\leq \epsilon_{i_0},\quad \nicefrac{\tilde{a}_{bi}}{a_{bi}}\leq \epsilon_{i_0},\quad \nicefrac{a_{iw}}{\tilde{a}_{iw}}\leq \epsilon_{i_0},\quad \nicefrac{\tilde{a}_{iw}}{a_{iw}}\leq \epsilon_{i_0},\quad \nicefrac{a_{bw}}{\tilde{a}_{bw}}\leq \epsilon_{i_0},\quad \nicefrac{\tilde{a}_{bw}}{a_{bw}}\leq\epsilon_{i_0}
		\end{equation*}
		for all $i\in D$, then it will imply $\eta^*\leq \epsilon_{i_0}$. Combining this with Proposition \ref{lower1}, we get $\eta^*=\epsilon_{i_0}$ and consequently, $(\tilde{A}_{bo},\tilde{A}_{ow})$ is an optimally modified PCS.\\\\
		It is clear that $\frac{\tilde{a}_{bw}}{a_{bw}}<\frac{a_{bw}}{\tilde{a}_{bw}}=\epsilon_{i_0}$. Consider $i\in D$. Then there are two possibilities.
		\begin{enumerate}[(i)]
			\item $a_{bi}\times a_{iw}\leq\frac{a_{bw}}{\epsilon_{i_0}}$.\\
			Observe that, in this case, we have $\frac{a_{bi}}{\tilde{a}_{bi}}\leq \frac{\tilde{a}_{bi}}{a_{bi}}=\bigg(\frac{a_{bw}}{\epsilon_{i_0}\times a_{bi}\times  a_{iw}}\bigg)^{\frac{1}{2}}$. Since $a_{bi_0}\times a_{i_0w}\leq a_{bi}\times a_{iw}$, we get $\frac{a_{bw}}{a_{bi}\times a_{iw}}\leq \frac{a_{bw}}{a_{bi_0}\times a_{i_0w}}=\epsilon_{i_0}^3$. This gives $\bigg(\frac{a_{bw}}{\epsilon_{i_0}\times a_{bi}\times  a_{iw}}\bigg)^{\frac{1}{2}}\leq \epsilon_{i_0}$, i.e., $\frac{a_{bi}}{\tilde{a}_{bi}}\leq \frac{\tilde{a}_{bi}}{a_{bi}}\leq \epsilon_{i_0}$. Similarly, we get $\frac{a_{iw}}{\tilde{a}_{iw}}\leq \frac{\tilde{a}_{iw}}{a_{iw}}\leq \epsilon_{i_0}$.
			\item $a_{bi}\times a_{iw}>\frac{a_{bw}}{\epsilon_{i_0}}.$\\
			First observe that $\max\{\epsilon_i,\epsilon_{i,j}:i,j\in D\}=\max \{\epsilon_{i_0},\epsilon_{j_0},\epsilon_{i_0,j_0}\}=\epsilon_{i_0}$. Also observe that, in this case, we have $\frac{\tilde{a}_{bi}}{a_{bi}}<\frac{a_{bi}}{\tilde{a}_{bi}}=\bigg(\frac{\epsilon_{i_0}\times a_{bi}\times  a_{iw}}{a_{bw}}\bigg)^{\frac{1}{2}}$. From \eqref{CV1}, we have $\frac{a_{bw}}{\epsilon_{i_0}} =\epsilon_{i_0}^2\times a_{bi_0}\times a_{i_0w}$. Now $\epsilon_{i,i_0}\leq \epsilon_{i_0}$ implies that $\frac{a_{bw}}{\epsilon_{i_0}} \geq \epsilon_{i,i_0}^2\times a_{bi_0}\times a_{i_0w}$. Using \eqref{CV2}, we get $\frac{a_{bw}}{\epsilon_{i_0}} \geq\frac{a_{bi}}{\epsilon_{i,i_0}}\times \frac{a_{iw}}{\epsilon_{i,i_0}}$. This gives $\frac{\tilde{a}_{bi}}{a_{bi}}<\frac{a_{bi}}{\tilde{a}_{bi}}=\bigg(\frac{\epsilon_{i_0}\times a_{bi}\times  a_{iw}}{a_{bw}}\bigg)^{\frac{1}{2}}\leq \epsilon_{i,i_0}\leq \epsilon_{i_0}$. Similarly, we can prove $\frac{\tilde{a}_{iw}}{a_{iw}}<\frac{a_{iw}}{\tilde{a}_{iw}}\leq \epsilon_{i_0}$. Hence the theorem.
		\end{enumerate}
	\end{proof}
	From Theorem \ref{exact_obj}, we get
	\begin{equation}\label{ana_obj}
		\begin{split}
			\epsilon^*&=\eta^*=\max \{\epsilon_{i_0},\epsilon_{j_0},\epsilon_{i_0,j_0}\}=\max\{\epsilon_i,\epsilon_{i,j}:i,j\in D\}\\
			&=\max\biggl\{\bigg(\frac{a_{bi}\times a_{iw}}{a_{bw}}\bigg)^{\frac{1}{3}},\bigg(\frac{a_{bw}}{a_{bi}\times a_{iw}}\bigg)^{\frac{1}{3}},\bigg(\frac{a_{bi}\times a_{iw}}{a_{bj}\times a_{jw}}\bigg)^{\frac{1}{4}},\bigg(\frac{a_{bj}\times a_{jw}}{a_{bi}\times a_{iw}}\bigg)^{\frac{1}{4}}:i,j\in D\biggr\},
		\end{split}
	\end{equation}
	which is an analytic form of $\epsilon^*$.
	\begin{proposition}\label{fix_values}
		Let $i_0,j_0$ be as above, and let $(\tilde{A}_{bo}^*,\tilde{A}_{ow}^*)$ be an optimally modified PCS. Then 
		\begin{equation}\label{fix}
			\begin{cases}
				\begin{cases}
					\tilde{a}_{bi_0}^*=\epsilon_{i_0}\times a_{bi_0}\\
					\tilde{a}_{i_0w}^*=\epsilon_{i_0}\times a_{i_0w}, \ \ \text{if } \epsilon_{i_0}=\max \{\epsilon_{i_0},\epsilon_{j_0},\epsilon_{i_0,j_0}\},\\
					\tilde{a}_{bw}^*=\frac{a_{bw}}{\epsilon_{i_0}},
				\end{cases}\\
				\begin{cases}
					\tilde{a}_{bj_0}^*=\frac{a_{bj_0}}{\epsilon_{j_0}}\\
					\tilde{a}_{j_0w}^*=\frac{a_{j_0w}}{\epsilon_{j_0}}, \ \ \quad \quad  \text{if } \epsilon_{j_0}=\max \{\epsilon_{i_0},\epsilon_{j_0},\epsilon_{i_0,j_0}\},\\
					\tilde{a}_{bw}^*=\epsilon_{j_0}\times a_{bw},
				\end{cases}\\
				\begin{cases}
					\tilde{a}_{bi_0}^*=\epsilon_{i_0,j_0}\times a_{bi_0}\\
					\tilde{a}_{i_0w}^*=\epsilon_{i_0,j_0}\times a_{i_0w}\\
					\tilde{a}_{bj_0}^*=\frac{a_{bj_0}}{\epsilon_{i_0,j_0}},\ \ \quad \quad \text{if } \epsilon_{i_0,j_0}=\max \{\epsilon_{i_0},\epsilon_{j_0},\epsilon_{i_0,j_0}\},\\
					\tilde{a}_{j_0w}^*=\frac{a_{j_0w}}{\epsilon_{i_0,j_0}}\\
					\tilde{a}_{bw}^*=(\epsilon_{i_0,j_0}\times a_{bi_0})\times (\epsilon_{i_0,j_0}\times a_{i_0w}).
				\end{cases}
			\end{cases}
		\end{equation}
	\end{proposition}
	\begin{proof}
		Since all three cases are similar, we shall discuss only one case. Let $\epsilon_{i_0}=\max \{\epsilon_{i_0},\epsilon_{j_0},\epsilon_{i_0,j_0}\}$. So, by Theorem \ref{exact_obj}, we have $\eta^*=\epsilon_{i_0}$. As discussed in Proposition \ref{lower1}, we get $1\leq\eta_{bi_0},\eta_{i_0w},\eta_{bw}\leq \eta^*=\epsilon_{i_0}$ such that $\tilde{a}_{bi_0}^*\in\{\frac{a_{bi_0}}{\eta_{bi_0}},\eta_{bi_0}\times a_{bi_0}\}$, $\tilde{a}_{i_0w}^*\in\{\frac{a_{i_0w}}{\eta_{i_0w}},\eta_{i_0w}\times a_{i_0w}\}$ and $\tilde{a}_{bw}^*\in\{\frac{a_{bw}}{\eta_{bw}},\eta_{bw}\times a_{bw}\}$. It is easy to verify that if one of $\eta_{bi_0},\eta_{i_0w}$ and $\eta_{bw}$ is strictly less that $\epsilon_{i_0}$, then at least one of other two values is strictly greater than $\epsilon_{i_0}$, which is not possible. So, $\eta_{bi_0}=\eta_{i_0w}=\eta_{bw}=\epsilon_{i_0}$. It is clear that among the eight possibilities of $(\tilde{a}_{bi_0}^*,\tilde{a}_{i_0w}^*,\tilde{a}_{bw}^*)$, $(\epsilon_{i_0}\times a_{bi_0},\epsilon_{i_0}\times a_{i_0w},\frac{a_{bw}}{\epsilon_{i_0}})$ is the only possibility for which $\tilde{a}_{bi_0}^*\times \tilde{a}_{i_0w}^*=\tilde{a}_{bw}^*$ holds. This completes the proof.
	\end{proof}
	Proposition \ref{fix_values} implies that, in any case, $\tilde{a}_{bw}^*$ remains same for all $(\tilde{A}_{bo}^*,\tilde{A}_{ow}^*)$.
	\begin{theorem}\label{opt_sol_PCS}
		If $(\tilde{A}_{bo}^*,\tilde{A}_{ow}^*)$ is an optimally modified PCS, then $\tilde{a}_{bw}^*$ is as defined by \eqref{fix}, $\tilde{a}_{bi}^*\in \bigg[\max\{\frac{a_{bi}}{\eta^*},\frac{\tilde{a}_{bw}^*}{\eta^*\times a_{iw}}\},\min\{\eta^*\times a_{bi},\frac{\eta^*\times \tilde{a}_{bw}^*}{a_{iw}}\}\bigg]$ and $\tilde{a}_{iw}^*=\frac{\tilde{a}_{bw}^*}{\tilde{a}_{bi}^*}$ (or equivalently, $\tilde{a}_{bi}^*=\frac{\tilde{a}_{bw}^*}{\tilde{a}_{iw}^*}$ and $\tilde{a}_{iw}^*\in \bigg[\max\{\frac{a_{iw}}{\eta^*},\frac{\tilde{a}_{bw}^*}{\eta^*\times a_{bi}}\},\min\{\eta^*\times a_{iw},\frac{\eta^*\times \tilde{a}_{bw}^*}{a_{bi}}\}\bigg]$) for all $i\in D$. Converse is also true.
	\end{theorem}
	\begin{proof}
		Let $(\tilde{A}_{bo}^*,\tilde{A}_{ow}^*)$ be an optimally modified PCS. Then, by Proposition \ref{fix_values}, $\tilde{a}_{bw}^*$ defined by \eqref{fix} is the only choice. Also, $\frac{a_{bi}}{\tilde{a}_{bi}^*},\frac{\tilde{a}_{bi}^*}{a_{bi}},\frac{a_{iw}}{\tilde{a}_{iw}^*},\frac{\tilde{a}_{iw}^*}{a_{iw}}\leq \eta^*$ for all $i\in D$. This gives 
		\begin{equation}\label{interval_1}
			\tilde{a}_{bi}^*\in \bigg[\frac{a_{bi}}{\eta^*},\eta^*\times a_{bi}\bigg] \text{ and } \tilde{a}_{iw}^*\in \bigg[\frac{a_{iw}}{\eta^*},\eta^*\times a_{iw}\bigg] \text{ for all }i\in D.
		\end{equation}
		Since $(\tilde{A}_{bo}^*,\tilde{A}_{ow}^*)$ is consistent, we get $\tilde{a}_{bi}^*=\frac{\tilde{a}_{bw}^*}{\tilde{a}_{iw}^*}$ and $\tilde{a}_{iw}^*=\frac{\tilde{a}_{bw}^*}{\tilde{a}_{bi}^*}$. This gives
		\begin{equation}\label{interval_2}
			\tilde{a}_{bi}^*\in \bigg[\frac{\tilde{a}_{bw}^*}{\eta^*\times a_{iw}},\frac{\eta^*\times \tilde{a}_{bw}^*}{a_{iw}}\bigg] \text{ and } \tilde{a}_{iw}^*\in \bigg[\frac{\tilde{a}_{bw}^*}{\eta^*\times a_{bi}},\frac{\eta^*\times \tilde{a}_{bw}^*}{a_{bi}}\bigg] \text{ for all } i\in D.
		\end{equation}
		Combining \eqref{interval_1} and \eqref{interval_2}, we get
		\begin{equation}\label{interval_3}
			\begin{split}
				&\tilde{a}_{bi}^*\in \bigg[\max\{\frac{a_{bi}}{\eta^*},\frac{\tilde{a}_{bw}^*}{\eta^*\times a_{iw}}\},\min\{\eta^*\times a_{bi},\frac{\eta^*\times \tilde{a}_{bw}^*}{a_{iw}}\}\bigg] \text{ and }\\
				&\tilde{a}_{iw}^*\in \bigg[\max\{\frac{a_{iw}}{\eta^*},\frac{\tilde{a}_{bw}^*}{\eta^*\times a_{bi}}\},\min\{\eta^*\times a_{iw},\frac{\eta^*\times \tilde{a}_{bw}^*}{a_{bi}}\}\bigg] \text{ for all }i\in D.
			\end{split}
		\end{equation} 
		Since we have $\tilde{a}_{bi}^*\times \tilde{a}_{iw}^*=\tilde{a}_{bw}^*$, we have choice for only one of $\tilde{a}_{bi}^*$ and $\tilde{a}_{iw}^*$. So, if we choose $\tilde{a}_{bi}^*$ satisfying \eqref{interval_3}, then $\tilde{a}_{iw}^*=\frac{\tilde{a}_{bw}^*}{\tilde{a}_{bi}^*}$ and vice versa.\\
		Conversely, let $(\tilde{A}_{bo},\tilde{A}_{ow})$ be a PCS such that the hypothesis is satisfied. Observe that $(\tilde{A}_{bo},\tilde{A}_{ow})$ is consistent. From the definition of $\tilde{a}_{bw}$ and Theorem \ref{exact_obj}, it is clear that $\frac{a_{bw}}{\tilde{a}_{bw}},\frac{\tilde{a}_{bw}}{a_{bw}}\leq\eta^*$. From the hypothesis, we get $\tilde{a}_{bi}\in \big[\frac{a_{bi}}{\eta^*},\eta^*\times a_{bi}\big]$ as well as $\tilde{a}_{bi}\in \big[\frac{\tilde{a}_{bw}}{\eta^*\times a_{iw}},\frac{\eta^*\times \tilde{a}_{bw}}{a_{iw}}\big]$, i.e.,  $\frac{\tilde{a}_{bw}}{\tilde{a}_{iw}}\in \big[\frac{\tilde{a}_{bw}}{\eta^*\times a_{iw}},\frac{\eta^*\times \tilde{a}_{bw}}{a_{iw}}\big]$ for all $i\in D$. This gives $\frac{a_{bi}}{\tilde{a}_{bi}},\frac{\tilde{a}_{bi}}{a_{bi}},\frac{a_{iw}}{\tilde{a}_{iw}},\frac{\tilde{a}_{iw}}{a_{iw}}\leq \eta^*$ for all $i\in D$. So, $(\tilde{A}_{bo},\tilde{A}_{ow})$ is an optimally modified PCS. Hence the proof.
	\end{proof}
	Theorem \ref{opt_sol_PCS} gives an analytic form of optimally modified PCS. Now, by \cite[Theorem 4]{wu2023analytical}, the collection of all optimal weights of criterion $c_i$ is $[{w_i^l}^*,{w_i^u}^*]$, where
	\begin{equation}\label{optimal_weights}
		\begin{split}
			{w_i^l}^*=\frac{\inf\{\tilde{a}_{iw}^*\}}{\inf\{\tilde{a}_{iw}^*\}+\displaystyle\sum_{j\neq i}\sup\{\tilde{a}_{jw}^*\}}&=\frac{\max\{\frac{a_{iw}}{\eta^*},\frac{\tilde{a}_{bw}^*}{\eta^*\times a_{bi}}\}}{\max\{\frac{a_{iw}}{\eta^*},\frac{\tilde{a}_{bw}^*}{\eta^*\times a_{bi}}\}+\displaystyle\sum_{j\neq i}\min\{\eta^*\times a_{jw},\frac{\eta^*\times \tilde{a}_{bw}^*}{a_{bj}}\}}\\
			&=\frac{\max\{\frac{a_{iw}}{\epsilon^*},\frac{\tilde{a}_{bw}^*}{\epsilon^*\times a_{bi}}\}}{\max\{\frac{a_{iw}}{\epsilon^*},\frac{\tilde{a}_{bw}^*}{\epsilon^*\times a_{bi}}\}+\displaystyle\sum_{j\neq i}\min\{\epsilon^*\times a_{jw},\frac{\epsilon^*\times \tilde{a}_{bw}^*}{a_{bj}}\}}\text{ and}\\
			{w_i^u}^*=\frac{\sup\{\tilde{a}_{iw}^*\}}{\sup\{\tilde{a}_{iw}^*\}+\displaystyle\sum_{j\neq i}\inf\{\tilde{a}_{jw}^*\}}&=\frac{\min\{\eta^*\times a_{iw},\frac{\eta^*\times \tilde{a}_{bw}^*}{a_{bi}}\}}{\min\{\eta^*\times a_{iw},\frac{\eta^*\times \tilde{a}_{bw}^*}{a_{bi}}\}+\displaystyle\sum_{j\neq i}\max\{\frac{a_{jw}}{\eta^*},\frac{\tilde{a}_{bw}^*}{\eta^*\times a_{bj}}\}}\\
			&=\frac{\min\{\epsilon^*\times a_{iw},\frac{\epsilon^*\times \tilde{a}_{bw}^*}{a_{bi}}\}}{\min\{\epsilon^*\times a_{iw},\frac{\epsilon^*\times \tilde{a}_{bw}^*}{a_{bi}}\}+\displaystyle\sum_{j\neq i}\max\{\frac{a_{jw}}{\epsilon^*},\frac{\tilde{a}_{bw}^*}{\epsilon^*\times a_{bj}}\}}, \ i\in C.
		\end{split}
	\end{equation}
	Equations \eqref{optimal_weights} and \eqref{ana_obj} give an analytic form of optimal solutions of problem \eqref{optimization_2}.
	\subsection{A secondary objective function to obtain unique optimal weight set}
	From equation \eqref{optimal_weights}, it follows that the non-linear multiplicative BWM may lead to multiple optimal weight sets. Although multiple weight sets are preferable in some cases, in most cases, the decision maker prefers a unique weight set. In this subsection, our goal is to select the best optimal weight set, for which we first obtain the best optimally modified PCS, which is done by introducing a secondary objective function.\\\\
	Inline with the core idea of \cite[subsection 3.3]{wu2023analytical}, an optimally modified PCS having the minimum distance between $(\tilde{a}_{bi}^*,\tilde{a}_{iw}^*)$ and $(\tilde{a}_{bi},\tilde{a}_{iw})$, i.e., the minimum value of $$\max\{\nicefrac{a_{bi}}{\tilde{a}_{bi}^*},\nicefrac{\tilde{a}_{bi}^*}{a_{bi}},\nicefrac{a_{iw}}{\tilde{a}_{iw}^*},\nicefrac{\tilde{a}_{iw}^*}{a_{iw}}\},$$ for all $i\in C$ is selected as the best optimally modified PCS.\\\\
	Consider the following optimization problem.
	\begin{equation}\label{optimization_7}
		\begin{split}
			&\min \max\{\nicefrac{a_{bj}}{\tilde{a}_{bj}},\nicefrac{\tilde{a}_{bj}}{a_{bj}},\nicefrac{a_{jw}}{\tilde{a}_{jw}},\nicefrac{\tilde{a}_{jw}}{a_{jw}}\}\\
			&\text{sub to:} \\
			&\nicefrac{a_{bi}}{\tilde{a}_{bi}}\leq \eta^*,\quad \nicefrac{\tilde{a}_{bi}}{a_{bi}}\leq \eta^*,\quad \nicefrac{a_{iw}}{\tilde{a}_{iw}}\leq \eta^*,\quad \nicefrac{\tilde{a}_{iw}}{a_{iw}}\leq \eta^*,\quad \nicefrac{a_{bw}}{\tilde{a}_{bw}}\leq \eta^*,\quad \nicefrac{\tilde{a}_{bw}}{a_{bw}}\leq\eta^*,\\
			&\tilde{a}_{bi}\times \tilde{a}_{iw}=\tilde{a}_{bw},\quad \tilde{a}_{bi},\tilde{a}_{iw},\tilde{a}_{bw}\geq 0 \text{ for all } i\in D.
		\end{split}
	\end{equation}
	Observe that problem \eqref{optimization_7} is a mini-max problem with $2n-3$ variables $\tilde{a}_{b1},...,\tilde{a}_{bb-1},\tilde{a}_{bb+1},...,\tilde{a}_{bn},$ $\tilde{a}_{1w},...,\tilde{a}_{w-1w},\tilde{a}_{w+1w},...,\tilde{a}_{nw}$ having the collection of all optimally modified PCS as the feasible region. So, it has optimal solution(s) of the form \\$(\tilde{a}_{b1}^*,...,\tilde{a}_{bb-1}^*,\tilde{a}_{bb+1}^*,...,\tilde{a}_{bn}^*,\tilde{a}_{1w}^*,...,\tilde{a}_{w-1w}^*,\tilde{a}_{w+1w}^*,...,\tilde{a}_{nw}^*)$. Each optimal solution along with $\tilde{a}_{bb}^*=\tilde{a}_{ww}^*=1$ is an optimally modified PCS having the minimum distance between $(\tilde{a}_{bj}^*,\tilde{a}_{jw}^*)$ and $(\tilde{a}_{bj},\tilde{a}_{jw})$ and this minimum value is precisely the optimal objective value.\\\\
	Note that problem \eqref{optimization_7} is equivalent to the following minimization problem. 
	\begin{equation}\label{optimization_8}
		\begin{split}
			&\min \eta_j\\
			&\text{sub to:} \\
			&\nicefrac{a_{bi}}{\tilde{a}_{bi}}\leq \eta^*,\quad \nicefrac{\tilde{a}_{bi}}{a_{bi}}\leq \eta^*,\quad \nicefrac{a_{iw}}{\tilde{a}_{iw}}\leq \eta^*,\quad \nicefrac{\tilde{a}_{iw}}{a_{iw}}\leq \eta^*,\quad \nicefrac{a_{bw}}{\tilde{a}_{bw}}\leq \eta^*,\quad \nicefrac{\tilde{a}_{bw}}{a_{bw}}\leq\eta^*,\\
			&\nicefrac{a_{bj}}{\tilde{a}_{bj}}\leq \eta_j,\quad \nicefrac{\tilde{a}_{bj}}{a_{bj}}\leq \eta_j,\quad \nicefrac{a_{iw}}{\tilde{a}_{iw}}\leq\eta_j,\quad \nicefrac{\tilde{a}_{iw}}{a_{iw}}\leq \eta_j,\\
			&\tilde{a}_{bi}\times \tilde{a}_{iw}=\tilde{a}_{bw},\quad \tilde{a}_{bi},\tilde{a}_{iw},\tilde{a}_{bw}\geq 0 \text{ for all } i\in D.
		\end{split}
	\end{equation}
	Here, the optimal objective value $\eta_j^*$ is the minimum possible distance between $(\tilde{a}_{bj}^*,\tilde{a}_{jw}^*)$ and $(\tilde{a}_{bj},\tilde{a}_{jw})$. Note that $\eta_j^*\leq \eta^*$. After calculating $\eta_j^*$ for all $j\in C$, the best optimally modified PCS is obtained using the following minimization problem.
	\begin{equation}\label{optimization_9}
		\begin{split}
			&\min f(=1)\\
			&\text{sub to:} \\
			&\nicefrac{a_{bj}}{\tilde{a}_{bj}}\leq \eta_j^*,\quad \nicefrac{\tilde{a}_{bj}}{a_{bj}}\leq \eta_j^*,\quad \nicefrac{a_{jw}}{\tilde{a}_{jw}}\leq\eta_j^*,\quad \nicefrac{\tilde{a}_{jw}}{a_{jw}}\leq \eta_j^*,\\
			&\tilde{a}_{bi}\times \tilde{a}_{iw}=\tilde{a}_{bw},\quad \tilde{a}_{bi},\tilde{a}_{iw},\tilde{a}_{bw}\geq 0 \text{ for all } i\in D,\ j\in C.
		\end{split}
	\end{equation}
	Note that problem \eqref{optimization_9} is a non-linear problem having $2n-3$ variables $\tilde{a}_{b1},...,\tilde{a}_{bb-1},\tilde{a}_{bb+1},...,\tilde{a}_{bn},$ $\tilde{a}_{1w},...,\tilde{a}_{w-1w},\tilde{a}_{w+1w},...,\tilde{a}_{nw}$. Now, Theorem \ref{unique} implies that the feasible region of this problem is a singleton set. So, the optimal solution\\ $(\tilde{a}_{b1}^*,...,\tilde{a}_{bb-1}^*,\tilde{a}_{bb+1}^*,...,\tilde{a}_{bn}^*,\tilde{a}_{1w}^*,...,\tilde{a}_{w-1w}^*,\tilde{a}_{w+1w}^*,...,\tilde{a}_{nw}^*)$ along with $\tilde{a}_{bb}^*=\tilde{a}_{ww}^*=1$ forms the best optimally modified PCS, which is unique. Observe that the objective function has no role in obtaining this PCS and thus any function can serve as an objective function. Here, we have taken the constant function $1$ as the objective function. 
	\begin{theorem}\label{unique}
		The following statements hold.
		\begin{enumerate}
			\item If $\eta^*=\epsilon_{i_0}$, then $(\tilde{A}_{bo}^*,\tilde{A}_{ow}^*)$ given by \eqref{m2} is the only PCS in the feasible region of problem \eqref{optimization_9}.
			\item If $\eta^*=\epsilon_{j_0}$, then $(\tilde{A}_{bo}^*,\tilde{A}_{ow}^*)$ given by \eqref{m1} is the only PCS in the feasible region of problem \eqref{optimization_9}.
			\item If $\eta^*=\epsilon_{i_0,j_0}$, then $(\tilde{A}_{bo}^*,\tilde{A}_{ow}^*)$ given by \eqref{m3} is the only PCS in the feasible region of problem \eqref{optimization_9}.
		\end{enumerate}
	\end{theorem}
	\begin{proof}
		Here, we shall prove only $1^{st}$ statement. Since $\eta_j^*\leq \eta^*$, any PCS other than optimally modified PCS does not belongs to the feasible region of problem \eqref{optimization_9}. Let $(\tilde{A}_{bo}'^*,\tilde{A}_{ow}'^*)$ be an optimally modified PCS. From Proposition \ref{fix_values}, we have $\tilde{a}_{bw}^*=\tilde{a}_{bw}'^*=\frac{a_{bw}}{\epsilon_{i_0}}$. So, $\max\{\nicefrac{a_{bw}}{\tilde{a}_{bw}^*},\nicefrac{\tilde{a}_{bw}^*}{a_{bw}}\}=\max\{\nicefrac{a_{bw}}{\tilde{a}_{bw}'^*},\nicefrac{\tilde{a}_{bw}'^*}{a_{bw}}\}$. Now, as discussed in Proposition \ref{lower1}, we get $1\leq\eta_{bi},\eta_{iw}\leq \epsilon_{i_0}$ such that $\max\{\nicefrac{a_{bi}}{\tilde{a}_{bi}'^*},\nicefrac{\tilde{a}_{bi}'^*}{a_{bi}}\}=\eta_{bi}$ and $\max\{\nicefrac{a_{iw}}{\tilde{a}_{iw}'^*,\nicefrac{\tilde{a}_{iw}'^*}{a_{iw}}\}=\eta_{iw}}$ for all $i\in D$. So, we get $\tilde{a}_{bi}'^*\in\{\frac{a_{bi}}{\eta_{bi}},\eta_{bi}\times a_{bi}\}$ and $\tilde{a}_{iw}'^*\in\{\frac{a_{iw}}{\eta_{iw}},\eta_{iw}\times a_{iw}\}$. Let $i\in D$ be such that $a_{bi}\times a_{iw}\leq \frac{a_{bw}}{\epsilon_{i_0}}$. Then we have $\frac{a_{bi}}{\tilde{a}_{bi}^*}\leq \frac{\tilde{a}_{bi}^*}{a_{bi}}=\bigg(\frac{a_{bw}}{\epsilon_{i_0}\times a_{bi}\times  a_{iw}}\bigg)^{\frac{1}{2}}$ and $\frac{a_{iw}}{\tilde{a}_{iw}^*}\leq \frac{\tilde{a}_{iw}^*}{a_{iw}}=\bigg(\frac{a_{bw}}{\epsilon_{i_0}\times a_{bi}\times  a_{iw}}\bigg)^{\frac{1}{2}}$. So, $\max\{\frac{a_{bi}}{\tilde{a}_{bi}^*},\frac{\tilde{a}_{bi}^*}{a_{bi}},\frac{a_{iw}}{\tilde{a}_{iw}^*},\frac{\tilde{a}_{iw}^*}{a_{iw}}\}=\bigg(\frac{a_{bw}}{\epsilon_{i_0}\times a_{bi}\times  a_{iw}}\bigg)^{\frac{1}{2}}$. Now, $\tilde{a}_{bi}'^*\times \tilde{a}_{iw}'^*=\frac{a_{bw}}{\epsilon_{i_0}}$ implies that if one of $\eta_{bi}$ and $\eta_{iw}$ is strictly less than $\bigg(\frac{a_{bw}}{\epsilon_{i_0}\times a_{bi}\times  a_{iw}}\bigg)^{\frac{1}{2}}$, then the other one is strictly greater than $\bigg(\frac{a_{bw}}{\epsilon_{i_0}\times a_{bi}\times  a_{iw}}\bigg)^{\frac{1}{2}}$. This gives $\max\{\frac{a_{bi}}{\tilde{a}_{bi}^*},\frac{\tilde{a}_{bi}^*}{a_{bi}},\frac{a_{iw}}{\tilde{a}_{iw}^*},\frac{\tilde{a}_{iw}^*}{a_{iw}}\}\leq \max\{\eta_{bi},\eta_{iw}\}=\max\{\frac{a_{bi}}{\tilde{a}_{bi}'^*},\frac{\tilde{a}_{bi}'^*}{a_{bi}},\frac{a_{iw}}{\tilde{a}_{iw}'^*},\frac{\tilde{a}_{iw}'^*}{a_{iw}}\}$, and if $\max\{\frac{a_{bi}}{\tilde{a}_{bi}^*},\frac{\tilde{a}_{bi}^*}{a_{bi}},\frac{a_{iw}}{\tilde{a}_{iw}^*},\frac{\tilde{a}_{iw}^*}{a_{iw}}\}=\max\{\frac{a_{bi}}{\tilde{a}_{bi}'^*},\frac{\tilde{a}_{bi}'^*}{a_{bi}},\frac{a_{iw}}{\tilde{a}_{iw}'^*},\frac{\tilde{a}_{iw}'^*}{a_{iw}}\}$, then $\eta_{bi}=\eta_{iw}=\bigg(\frac{a_{bw}}{\epsilon_{i_0}\times a_{bi}\times  a_{iw}}\bigg)^{\frac{1}{2}}$, which gives $\tilde{a}_{bi}'^*=\tilde{a}_{bi}^*$ and $\tilde{a}_{iw}'^*=\tilde{a}_{iw}^*$. The same conclusion can be drawn if $a_{bi}\times a_{iw}<\frac{\tilde{a}_{bw}}{\epsilon_{i_0}}$. So, $\max\{\frac{a_{bi}}{\tilde{a}_{bi}^*},\frac{\tilde{a}_{bi}^*}{a_{bi}},\frac{a_{iw}}{\tilde{a}_{iw}^*},\frac{\tilde{a}_{iw}^*}{a_{iw}}\}\leq\max\{\frac{a_{bi}}{\tilde{a}_{bi}'^*},\frac{\tilde{a}_{bi}'^*}{a_{bi}},\frac{a_{iw}}{\tilde{a}_{iw}'^*},\frac{\tilde{a}_{iw}'^*}{a_{iw}}\}$ for all $i\in C$, and $(\tilde{A}_{bo}^*,\tilde{A}_{ow}^*)$ is the only PCS having this property. This completes the proof.
	\end{proof}
	Using the best optimally modified PCS, the best optimal weight set is calculated by \eqref{weights}, which is the resultant weight set.
	\subsection{Consistency analysis}
	The outcome of an MCDM method relies on the decision data which is usually irrational due to human involvement. Estimation of this discrepancy is an important aspect of any MCDM method. For BWM, inconsistency in preferences (and consequently, in resultant weights) is measured by a ratio called Consistency Ratio (CR) defined as
	\begin{equation}\label{CR}
		\text{CR}=\frac{\epsilon^*}{\text{Consistency Index (CI)}},
	\end{equation}
	where CI$=\sup\{\epsilon^*:\epsilon^*$ is the optimal objective value of problem \eqref{optimization_2} for some $(A_{bo},A_{ow})$ having the given value of $a_{bw}$\}\cite{rezaei2015best}. So, CI can be thought as a function of $a_{bw}$. In the existing literature, CR is also referred as inconsistency level\cite{brunelli2019multiplicative}. In this subsection, our goal is to obtain analytic forms of CI and CR.\\\\
	From equation \eqref{ana_obj}, it follows that to obtain CI$_{a_{bw}}$, it is sufficient to find the maximum possible values of $(\frac{a_{bi}\times a_{iw}}{a_{bw}})^{\frac{1}{3}},(\frac{a_{bw}}{a_{bi}\times a_{iw}})^{\frac{1}{3}},(\frac{a_{bi}\times a_{iw}}{a_{bj}\times a_{jw}})^{\frac{1}{4}}$ and $(\frac{a_{bj}\times a_{jw}}{a_{bi}\times a_{iw}})^{\frac{1}{4}}$, where $i,j\in D$. Now, observe that if $a_{bi_1}\leq a_{bi_2}$ and $a_{i_1w}\leq a_{i_2w}$ for some $i_1,i_2\in D$, then $(\frac{a_{bi_1}\times a_{i_1w}}{a_{bw}})^{\frac{1}{3}}\leq (\frac{a_{bi_2}\times a_{i_2w}}{a_{bw}})^{\frac{1}{3}}$ and $(\frac{a_{bw}}{a_{bi_1}\times a_{i_1w}})^{\frac{1}{3}}\geq (\frac{a_{bw}}{a_{bi_2}\times a_{i_2w}})^{\frac{1}{3}}$. So, $(\frac{a_{bi}\times a_{iw}}{a_{bw}})^{\frac{1}{3}}$ and $(\frac{a_{bw}}{a_{bi}\times a_{iw}})^{\frac{1}{3}}$ have the maximum possible values if $a_{bi}=a_{iw}=a_{bw}$ and $a_{bi}=a_{iw}=1$ for some $i\in D$, respectively. This gives that the maximum possible value of $(\frac{a_{bi}\times a_{iw}}{a_{bw}})^{\frac{1}{3}}$ and $(\frac{a_{bw}}{a_{bi}\times a_{iw}})^{\frac{1}{3}}$ is $a_{bw}^{\frac{1}{3}}$. Similarly, the maximum possible value of $(\frac{a_{bi}\times a_{iw}}{a_{bj}\times a_{jw}})^{\frac{1}{4}}$ and $(\frac{a_{bj}\times a_{jw}}{a_{bi}\times a_{iw}})^{\frac{1}{4}}$ is $a_{bw}^{\frac{1}{2}}$. So, we get CI$_{a_{bw}}=\max\{a_{bw}^{\frac{1}{3}},a_{bw}^{\frac{1}{2}}\}$. For $\frac{1}{9}$ to $9$ scale, we have $a_{bw}\geq 1$. So, $a_{bw}^\frac{1}{3}\leq a_{bw}^\frac{1}{2}$, which gives CI$_{a_{bw}}=a_{bw}^\frac{1}{2}$. This along with equation \eqref{ana_obj} gives
	\begin{equation}\label{ana_cr}
		\text{CR}=\frac{\max\biggl\{\bigg(\frac{a_{bi}\times a_{iw}}{a_{bw}}\bigg)^{\frac{1}{3}},\bigg(\frac{a_{bw}}{a_{bi}\times a_{iw}}\bigg)^{\frac{1}{3}},\bigg(\frac{a_{bi}\times a_{iw}}{a_{bj}\times a_{jw}}\bigg)^{\frac{1}{4}},\bigg(\frac{a_{bj}\times a_{jw}}{a_{bi}\times a_{iw}}\bigg)^{\frac{1}{4}}:i,j\in D\biggr\}}{a_{bw}^{\frac{1}{2}}}.
	\end{equation}
	Now, we discuss some key features of analytic forms of CI and CR.
	\begin{enumerate}
		\item As discussed in Subsection 2.2, multiplicative BWM is linearized using the logarithmic transformation. It is important to mention that the values of $ln(\text{CI})$ given in Table \ref{ia} matches perfectly with the CI of linearized model given in \cite[Table 1]{brunelli2019multiplicative}.
		\item Liang et al.\cite{liang2020consistency} has discussed the importance of consistency indicator which can provide immediate feedback to the decision maker about the irrationality of decision judgments. The output-based formulation of CR defined by \eqref{CR} lacks this ability. This limitation is resolved by the analytic form of CR defined by \eqref{ana_cr}, which is an input-based formulation of CR.
		\item It is established that for a consistency indicator to exhibit reasonable behavior, it must satisfy some specific properties\cite{brunelli2015axiomatic,koczkodaj2018axiomatization}. Some of these properties, outlined in Proposition \ref{cr_pro}, are stated without proof, as their proofs align with those of \cite[Proposition 1]{liang2020consistency}.
	\end{enumerate}
	\begin{proposition}\label{cr_pro}
		The CR satisfies the following properties.
		\begin{enumerate}
			\item CR$=0$ if and only if $(A_{bo},A_{ow})$ is consistent.
			\item CR is normalized, i.e., $0\leq$ CR $\leq 1$.
			\item CR is invariant with respect to a permutation of the indices of the criteria.
			\item If we remove $i^{th}$ criterion, where $i\neq b,w,i_0,j_0$, then the value of CR does not change.
			\item CR is a continuous function of $a_{bi},a_{iw}$ and $a_{bw}$.
			\item For a consistent $(A_{bo},A_{ow})$, moving one of the preferences $a_{bi}$ or $a_{iw}$ away from the original value in the range $[1,a_{bw}]$ will lead to an increase in the value of CR.
		\end{enumerate}
	\end{proposition}
	\begin{table}
		\caption{The values of CI and $ln(\text{CI})$}\label{ia}
		\begin{tabular}{c c c c c c c c c c}
			\hline
			$a_{bw}$&$1$&$2$&$3$&$4$&$5$&$6$&$7$&$8$&$9$\\
			\hline
			CI&$1$&$1.4142$&$1.7320$&$2$&$2.2361$&$2.4494$&$2.6457$&$2.8284$&$3$\\
			\hline
			$\ln(\text{CI})$&$0$&$0.3466$&$0.5493$&$0.6931$&$0.8047$&$0.8959$&$0.9729$&$1.0397$&$1.0986$\\
			\hline
		\end{tabular}
	\end{table}
	\subsection{Numerical Examples}
	In this subsection, we illustrate the proposed approach and compare it with the existing approaches using numerical examples. For $j\in C$, define $\eta_j=\max\{\nicefrac{a_{bj}}{\frac{w_b}{w_j}},\nicefrac{\frac{w_b}{w_j}}{a_{bj}},\nicefrac{a_{jw}}{\frac{w_j}{w_w}},\nicefrac{\frac{w_j}{w_w}}{a_{jw}}\}$, where $W=\{w_1,w_2,...,w_n\}$ is any weight set. It is clear that a lower value of $\eta_j$ indicates a better retention of the original decision data, and thus a better weight set.\\\\
	\textbf{Example 1: }Let $C=\{c_1,c_2,...,c_5\}$ be the set of decision criteria with $c_2$ as the best and $c_5$ as the worst criterion. Let $A_{bo}=(2,1,5,3,8)$ and $A_{ow}=(4,8,3,3,1)$ be the best-to-other and the other-to-worst vector, respectively.\\\\
	Here, $D_1=\phi$ and $j_0=3$. So, $\epsilon^*=\epsilon_3=1.2331$. By Proposition \ref{fix_values}, we get $\tilde{a}_{bw}^*=9.8648$. Now, by equation \eqref{optimal_weights}, optimal interval-weights are $w_1^*=[0.1905,0.2360]$, $w_2^*=[0.4498,0.4941]$, $w_3^*=[0.1109,0.1219]$, $w_4^*=[0.1276,0.1762]$ and $w_5^*=[0.0456,0.0501]$. Now, Theorem \ref{unique} implies that the PCS given by equation \eqref{m1} is the best optimally modified PCS. So, $\tilde{A}_{bo}^*=(2.2209,1,4.0548,3.1408,9.8648)$ and $\tilde{A}_{ow}^*=(4.4418,9.8648,2.4329,3.1408,1)$ forms the best optimally modified PCS and thus, by equation \eqref{weights}, $\{0.2127,0.4724,0.1165,0.1504,0.0479\}$ is the best optimal weight set. The comparison between this weight set and the centre of the interval-weight is given in Table \ref{example_1}.
	\begin{table}[H]
		\centering
		\caption{Computed weights: Example 1\label{example_1}}
		\begin{tabular}{||c||c|c|c||c|c||}
			\hline
			\hline
			\multirow{3}{*}{Criterion}&\multicolumn{3}{c||}{Non-linear multiplicative}&\multicolumn{2}{c||}{The proposed}\\
			&\multicolumn{3}{c||}{model}&\multicolumn{2}{c||}{model}\\
			\cline{2-6}
			&Interval-weight&Centre&$\eta_j$&Weight&$\eta_j$\\
			\hline
			\hline
			$c_1$&$[0.1905,0.2360]$&$0.2133$&$1.1144$&$0.2127$&$1.1104$\\
			$c_2$&$[0.4498,0.4941]$&$0.4720$&$1.2331$&$0.4724$&$1.2331$\\
			$c_3$&$[0.1109,0.1219]$&$0.1162$&$1.2331$&$0.1165$&$1.2331$\\
			$c_4$&$[0.1276,0.1762]$&$0.1519$&$1.0582$&$0.1504$&$1.0469$\\
			$c_5$&$[0.0456,0.0501]$&$0.0478$&$1.2331$&$0.0479$&$1.2331$\\
			\hline
			\hline				
		\end{tabular}
	\end{table}
	\textbf{Example 2: }Let $C=\{c_1,c_2,...,c_5\}$ be the set of decision criteria with $c_2$ as the best and $c_5$ as the worst criterion. Let $A_{bo}=(2,1,4,5,9)$ and $A_{ow}=(3,9,2,2,1)$ be the best-to-other and the other-to-worst vector, respectively.\\\\
	Here, $i_0=1$ and $j_0=4$. So, $\epsilon^*=\max\{\epsilon_1,\epsilon_4,\epsilon_{1,4}\}=\max\{1.1447,1.0357,1.1362\}=1.1447=\epsilon_1$. By Proposition \ref{fix_values}, we get $\tilde{a}_{bw}^*=7.8622$. Now, by equation \eqref{optimal_weights}, optimal interval-weights are $w_1^*=[0.2101,0.2175]$, $w_2^*=[0.4810,0.4979]$, $w_3^*=[0.1103,0.1381]$, $w_4^*=[0.1072,0.1136]$ and $w_5^*=[0.0612,0.0633]$. Now, Theorem \ref{unique} implies that the PCS given by equation \eqref{m2} is the best optimally modified PCS. So, $\tilde{A}_{bo}^*=(2.2894,1,3.9654,1.0357,7.8622)$ and $\tilde{A}_{ow}^*=(3.4341,7.8622,1.9827,1.7734,1)$ forms the best optimally modified PCS and thus, by equation \eqref{weights}, $\{0.2139,0.4898,0.1235,0.1105,0.0623\}$ is the best optimal weight set. The comparison between this weight set and the centre of the interval-weight is given in Table \ref{example_2}.
	\begin{table}[H]
		\centering
		\caption{Computed weights: Example 2\label{example_2}}
		\begin{tabular}{||c||c|c|c||c|c||}
			\hline
			\hline
			\multirow{3}{*}{Criterion}&\multicolumn{3}{c||}{Non-linear multiplicative}&\multicolumn{2}{c||}{The proposed}\\
			&\multicolumn{3}{c||}{model\cite{brunelli2019multiplicative}}&\multicolumn{2}{c||}{model}\\
			\cline{2-6}
			&Interval-weight&Centre&$\eta_j$&Weight&$\eta_j$\\
			\hline
			\hline
			$c_1$&$[0.2101,0.2175]$&$0.2138$&$1.1447$&$0.2139$&$1.1447$\\
			$c_2$&$[0.4810,0.4979]$&$0.4894$&$1.1447$&$0.4898$&$1.1447$\\
			$c_3$&$[0.1103,0.1381]$&$0.1242$&$1.0149$&$0.1235$&$1.0087$\\
			$c_4$&$[0.1072,0.1136]$&$0.1104$&$1.1280$&$0.1105$&$1.1278$\\
			$c_5$&$[0.0612,0.0633]$&$0.0622$&$1.1447$&$0.0623$&$1.1447$\\
			\hline
			\hline				
		\end{tabular}
	\end{table}
	\textbf{Example 3: }Let $C=\{c_1,c_2,c_3,c_4\}$ be the set of decision criteria with $c_2$ as the best and $c_4$ as the worst criterion. Let $A_{bo}=(1,1,3,4)$ and $A_{ow}=(2,4,4,1)$ be the best-to-other and the other-to-worst vector, respectively.\\\\
	Here, $i_0=1$ and $j_0=3$. So, $\epsilon^*=\max\{\epsilon_1,\epsilon_3,\epsilon_{1,3}\}=\max\{1.2599,1.4422,1.5651\}=1.5651=\epsilon_{1,3}$. By Proposition \ref{fix_values}, we get $\tilde{a}_{bw}^*=4.8990$. Now, by equation \eqref{optimal_weights}, interval-weights are $w_1^*=[0.2701,0.2701]$, $w_2^*=[0.4228,0.4228]$, $w_3^*=[0.2206,0.2206]$ and $w_4^*=[0.0863,0.0863]$. Now, Theorem \ref{unique} implies that the PCS given by equation \eqref{m3} is the best optimally modified PCS. So, $\tilde{A}_{bo}^*=(1.5651,1,1.9168,4.8990)$ and $\tilde{A}_{ow}^*=(3.1302,4.8990,2.5558,1)$ forms the best optimally modified PCS and thus, by equation \eqref{weights}, $\{0.2701,0.4228,0.2206,0.0863\}$ is the best optimal weight set. The comparison between this weight set and the centre of the interval-weight is given in Table \ref{example_3}.
	\begin{table}[H]
		\centering
		\caption{Computed weights: Example 3\label{example_3}}
		\begin{tabular}{||c||c|c|c||c|c||}
			\hline
			\hline
			\multirow{3}{*}{Criterion}&\multicolumn{3}{c||}{Non-linear multiplicative}&\multicolumn{2}{c||}{The proposed}\\
			&\multicolumn{3}{c||}{model\cite{brunelli2019multiplicative}}&\multicolumn{2}{c||}{model}\\
			\cline{2-6}
			&Interval-weight&Centre&$\eta_j$&Weight&$\eta_j$\\
			\hline
			\hline
			$c_1$&$[0.2701,0.2701]$&$0.2701$&$1.5651$&$0.2701$&$1.5651$\\
			$c_2$&$[0.4228,0.4228]$&$0.4228$&$1.2247$&$0.4228$&$1.2247$\\
			$c_3$&$[0.2206,0.2206]$&$0.2206$&$1.5651$&$0.2206$&$1.5651$\\
			$c_4$&$[0.0863,0.0863]$&$0.0863$&$1.2247$&$0.0863$&$1.2247$\\
			\hline
			\hline				
		\end{tabular}
	\end{table}
	In Example 1 and 2, the proposed model has lower values of $\eta_j$ than the centre of interval-weight for all for all $j\in C$, which indicates the superiority of the proposed model over the centre of the interval-weight. In Example 3, the non-linear multiplicative model has a unique solution, and thus, both approaches coincide.
	\section{A real-world application}
	In this section, we discuss a real-world application of the proposed model in ranking the drivers of sustainable additive manufacturing.\\\\
	Additive Manufacturing (AM) is the process of creating a physical object by adding material layer-by-layer, and is therefore also known as layered manufacturing. The process of AM usually begins by creating a digital model of an object using Computer Aided Design (CAD) or by scanning an existing object. This model is cut into horizontal layers, and then the manufacturing equipment creates the final product layer by layer using materials such as plastic polymers, metals, ceramic, wood, etc\cite{thompson2016design}. The AM offers several significant advantages over traditional manufacturing methods such as improved design flexibility\cite{alogla2021impact}, rapid prototyping\cite{khorasani2022additive}, on-demand production\cite{ahmed2020demand}, etc., and thus, AM is used in various fields like medical\cite{salmi2021additive}, food\cite{yu2023future}, aerospace\cite{khorasani2022additive}, defense\cite{salunkhe2023current}, robotics\cite{stano2021additive}, and so on. \\\\
	The practice of using AM in a way that maximizes resource efficiency and minimizes environmental impact is called Sustainable AM (SAM). Some of the key strategies related to SAM are selection of environmentally friendly materials\cite{al2020additive}, waste reduction\cite{javaid2021role}, localized production\cite{verboeket2021additive}, life cycle assessment\cite{garcia2021comparative}, etc. These strategies aim to promote sustainability throughout the manufacturing process. From an economic standpoint, SAM offers cost-saving benefits by consuming less material and energy compared to conventional manufacturing methods\cite{peng2018sustainability}. This reduction in operational costs is attributed to the optimized use of resources and streamlined production processes inherent in SAM practices. Additionally, SAM has the potential to streamline supply chains by reducing reliance on multiple suppliers, lowering transportation expenses, and simplifying inventory management, thereby improving overall operational efficiency and cost-effectiveness\cite{javaid2021role}. From an environmental perspective, a study have demonstrated a notable reduction in carbon dioxide emissions compared to products made through conventional manufacturing techniques\cite{javaid2021role}. Moreover, SAM processes such as Carbon Capture and Utilization (CCU) technologies have emerged as promising solutions for further reducing carbon emissions by capturing and repurposing carbon dioxide emissions generated during industrial processes\cite{bara2021carbon}. Furthermore, SAM tends to produce minimal waste, contrasting starkly with the substantial waste generated by conventional manufacturing techniques, thus contributing positively to environmental sustainability efforts\cite{javaid2021role}. On a societal level, SAM has the potential to revolutionize education by democratizing access to prototyping and experimentation\cite{alabi2019applications}. By making AM technologies more accessible and affordable, SAM enables students, educators, and researchers to engage in hands-on learning experiences, fostering creativity, innovation, and interdisciplinary collaboration within educational institutions and research facilities.\\\\
	To enable the successful implementation of SAM, it is essential to analyze the factors that drive SAM adoption. Some of these drivers are reusing and recycling materials and packaging\cite{MOKTADIR20181366}, collaborative decision making\cite{GOVINDAN2015182}, green innovation\cite{su8080824}. Agrawal and Vinodh\cite{agrawal2021prioritisation} identified 40 such drivers, which are divided into 8 categories as given in Table \ref{driver}. To weight these categories and the drivers from the same category, four experts are considered to provide the pairwise comparisons. These values are adopted from \cite{agrawal2021prioritisation}. The average weights of these drivers and their ranking is given in Table \ref{real}. The results of comparisons are the same as numerical example. The resultant weight set of the proposed approach has lower values of $\eta_j$ for all categories and drivers than the centre of interval-weight. Although it is important to mention that, for this particular application, the difference between these weight sets is very less.
	\begin{table}
		\centering
		\caption{Drivers of sustainable additive manufacturing\cite{agrawal2021prioritisation}\label{driver}}
		\begin{tabular}{|c|c|}
			\hline
			Category& Driver\\
			\hline
			\hline
			\multirow{3}{*}{Material ($c_1$)} & Optimum resources utilisation ($c_{11}$)\\
			&Reusing and recycling materials and packaging ($c_{12}$)\\
			&Less material wastage ($c_{13}$)\\
			\hline
			\multirow{5}{*}{Management and stake holders ($c_2$)} & Proper training and education ($c_{21}$)\\
			&Incentives ($c_{22}$)\\
			&Stake holder participation ($c_{23}$)\\
			&Management support ($c_{24}$)\\
			&Employee commitment towards sustainability ($c_{25}$) \\
			\hline
			\multirow{6}{*}{Customer/Supplier/Competitor ($c_3$)} & Customers' expectations ($c_{31}$)\\
			&Supplier commitment ($c_{32}$)\\
			&Company image ($c_{33}$)\\
			&Public pressure ($c_{34}$)\\
			&Supply chain pressure ($c_{35}$)\\
			&Competitor pressure towards greening ($c_{36}$)\\
			\hline
			\multirow{5}{*}{Collaboration and trends ($c_4$)} & Market trends ($c_{41}$)\\
			&Economic benefits ($c_{42}$)\\
			&Collaboration between organizations ($c_{43}$) \\
			&Long-term survival in market ($c_{44}$)\\
			&Collaborative decision making ($c_{45}$)\\
			\hline
			\multirow{6}{*}{Technology ($c_5$)} &Green innovation ($c_{51}$)\\
			&Energy conservation ($c_{52}$)\\
			&Technology competence ($c_{53}$)\\
			&Technology adaptability ($c_{54}$)\\
			&Time to develop new product ($c_{55}$)\\
			& Effective visual control ($c_{56}$)\\
			\hline
			\multirow{5}{*}{Standards and regulations ($c_6$)} &Environmental conservation ($c_{61}$)\\
			&Emissions and global climate ($c_{62}$)\\
			&Safety and security ($c_{63}$)\\
			&Operational safety ($c_{64}$)\\
			&Compliance with regulations ($c_{65}$)\\
			\hline
			\multirow{4}{*}{Design ($c_7$)} &Light weight parts ($c_{71}$)\\
			&Eco design ($c_{72}$) \\
			&Flexibility in producing geometry ($c_{73}$)\\
			&Ease of producing complex structure ($c_{74}$)\\
			\hline
			\multirow{6}{*}{Process ($c_8$)} &Cost saving ($c_{81}$)\\
			&Time saving ($c_{82}$)\\
			&On-demand manufacturing ($c_{83}$)\\
			&Enables mass customization ($c_{84}$)\\
			&No specialized tooling required ($c_{85}$)\\
			&Quality ($c_{86}$)\\
			\hline
		\end{tabular}
	\end{table}
	\begin{sidewaystable}[h]
		\caption {Computed weights of the drivers of SAM\label{real}}
		\centering
		\begin{adjustbox}{width=\textwidth}
			\small
			\begin{tabular}{|c|c|c|c|c|c|c|c|c|c|c|c|}
				\hline
				\multicolumn{6}{|c|}{Non-linear multiplicative model}&\multicolumn{6}{c|}{The proposed model}\\
				\hline
				Category&Centre of interval-weight&Driver&Local weight&Global weight&Ranking&Category&Weight&Driver&Local weight&Global weight&Ranking\\
				\hline
				\multirow{3}{*}{$c_1$}&\multirow{3}{*}{$0.1391$}&$c_{11}$&$0.4562$&$0.0635$&$3$&\multirow{3}{*}{$c_1$}&\multirow{3}{*}{$0.1391$}&$c_{11}$&$0.4562$&$0.0635$&$3$\\
				&&$c_{12}$&$0.3100$&$0.0431$&$7$&&&$c_{12}$&$0.3100$&$0.0431$&$7$\\
				&&$c_{13}$&$0.2337$&$0.0325$&$13$&&&$c_{13}$&$0.2337$&$0.0325$&$13$\\
				\hline
				\multirow{5}{*}{$c_2$}&\multirow{5}{*}{$0.0511$}&$c_{21}$&$0.2232$&$0.0114$&$30$&\multirow{5}{*}{$c_2$}&\multirow{5}{*}{$0.0511$}&$c_{21}$&$0.2232$&$0.0114$&$30$\\
				&&$c_{22}$&$0.0652$&$0.0033$&$40$&&&$c_{22}$&$0.0652$&$0.0033$&$40$\\
				&&$c_{23}$&$0.1713$&$0.0088$&$33$&&&$c_{23}$&$0.1713$&$0.0088$&$33$\\
				&&$c_{24}$&$0.2610$&$0.0133$&$27$&&&$c_{24}$&$0.2610$&$0.0133$&$27$\\
				&&$c_{25}$&$0.2793$&$0.0143$&$25$&&&$c_{25}$&$0.2793$&$0.0143$&$25$\\
				\hline
				\multirow{6}{*}{$c_3$}&\multirow{6}{*}{$0.0618$}&$c_{31}$&$0.2337$&$0.0144$&$24$&\multirow{6}{*}{$c_3$}&\multirow{6}{*}{$0.0618$}&$c_{31}$&$0.2337$&$0.0144$&$24$\\
				&&$c_{32}$&$0.1140$&$0.0070$&$36$&&&$c_{32}$&$0.1140$&$0.0070$&$36$\\
				&&$c_{33}$&$0.2041$&$0.0126$&$29$&&&$c_{33}$&$0.2041$&$0.0126$&$29$\\
				&&$c_{34}$&$0.1347$&$0.0083$&$35$&&&$c_{34}$&$0.1347$&$0.0083$&$35$\\
				&&$c_{35}$&$0.0826$&$0.0051$&$38$&&&$c_{35}$&$0.0826$&$0.0051$&$38$\\
				&&$c_{36}$&$0.2308$&$0.0143$&$26$&&&$c_{36}$&$0.2308$&$0.0143$&$26$\\
				\hline
				\multirow{5}{*}{$c_4$}&\multirow{5}{*}{$0.0834$}&$c_{41}$&$0.1023$&$0.0085$&$34$&\multirow{5}{*}{$c_4$}&\multirow{5}{*}{$0.0833$}&$c_{41}$&$0.1023$&$0.0085$&$34$\\
				&&$c_{42}$&$0.3206$&$0.0267$&$15$&&&$c_{42}$&$0.3206$&$0.0267$&$15$\\
				&&$c_{43}$&$0.2428$&$0.0202$&$17$&&&$c_{43}$&$0.2428$&$0.0202$&$17$\\
				&&$c_{44}$&$0.1299$&$0.0108$&$31$&&&$c_{44}$&$0.1299$&$0.0108$&$31$\\
				&&$c_{45}$&$0.2044$&$0.0170$&$22$&&&$c_{45}$&$0.2044$&$0.0170$&$22$\\
				\hline
				\multirow{6}{*}{$c_5$}&\multirow{6}{*}{$0.2250$}&$c_{51}$&$0.3020$&$0.0680$&$2$&\multirow{6}{*}{$c_5$}&\multirow{6}{*}{$0.2250$}&$c_{51}$&$0.3020$&$0.0680$&$2$\\
				&&$c_{52}$&$0.2394$&$0.0539$&$4$&&&$c_{52}$&$0.2394$&$0.0539$&$4$\\
				&&$c_{53}$&$0.1992$&$0.0448$&$6$&&&$c_{53}$&$0.1992$&$0.0448$&$6$\\
				&&$c_{54}$&$0.1192$&$0.0268$&$14$&&&$c_{54}$&$0.1192$&$0.0268$&$14$\\
				&&$c_{55}$&$0.0829$&$0.0187$&$18$&&&$c_{55}$&$0.0829$&$0.0187$&$18$\\
				&&$c_{56}$&$0.0572$&$0.0129$&$28$&&&$c_{56}$&$0.0572$&$0.0129$&$28$\\
				\hline
				\multirow{5}{*}{$c_6$}&\multirow{5}{*}{$0.0718$}&$c_{61}$&$0.3447$&$0.0247$&$16$&\multirow{5}{*}{$c_6$}&\multirow{5}{*}{$0.0718$}&$c_{61}$&$0.3447$&$0.0247$&$16$\\
				&&$c_{62}$&$0.2573$&$0.0184$&$20$&&&$c_{62}$&$0.2573$&$0.0184$&$20$\\
				&&$c_{63}$&$0.0721$&$0.0052$&$37$&&&$c_{63}$&$0.0721$&$0.0052$&$37$\\
				&&$c_{64}$&$0.0679$&$0.0049$&$39$&&&$c_{64}$&$0.0679$&$0.0049$&$39$\\
				&&$c_{65}$&$0.2580$&$0.0185$&$19$&&&$c_{65}$&$0.2580$&$0.0185$&$19$\\
				\hline
				\multirow{4}{*}{$c_7$}&\multirow{4}{*}{$0.1839$}&$c_{71}$&$0.2868$&$0.0527$&$5$&\multirow{4}{*}{$c_7$}&\multirow{4}{*}{$0.1839$}&$c_{71}$&$0.2868$&$0.0527$&$5$\\
				&&$c_{72}$&$0.4204$&$0.0773$&$1$&&&$c_{72}$&$0.4204$&$0.0773$&$1$\\
				&&$c_{73}$&$0.1978$&$0.0364$&$12$&&&$c_{73}$&$0.1978$&$0.0364$&$12$\\
				&&$c_{74}$&$0.0950$&$0.0175$&$21$&&&$c_{74}$&$0.0950$&$0.0175$&$21$\\
				\hline
				\multirow{6}{*}{$c_8$}&\multirow{6}{*}{$0.1839$}&$c_{81}$&$0.0802$&$0.0148$&$23$&\multirow{6}{*}{$c_8$}&\multirow{6}{*}{$0.1839$}&$c_{81}$&$0.0802$&$0.0148$&$23$\\
				&&$c_{82}$&$0.2080$&$0.0383$&$11$&&&$c_{82}$&$0.2080$&$0.0383$&$11$\\
				&&$c_{83}$&$0.2315$&$0.0426$&$8$&&&$c_{83}$&$0.2315$&$0.0426$&$8$\\
				&&$c_{84}$&$0.2171$&$0.0399$&$9$&&&$c_{84}$&$0.2171$&$0.0399$&$9$\\
				&&$c_{85}$&$0.2086$&$0.0384$&$10$&&&$c_{85}$&$0.2086$&$0.0384$&$10$\\
				&&$c_{86}$&$0.0546$&$0.0100$&$32$&&&$c_{86}$&$0.0546$&$0.0100$&$32$\\
				\hline				
			\end{tabular}	
		\end{adjustbox}
	\end{sidewaystable}
	\section{Conclusions and future direction}
	The BWM is one of the recent MCDM methods that has been widely applied in various real-life applications. By incorporating different distance functions with the core idea of the original BWM, several models of BWM are developed, of which the non-linear multiplicative model is the focus of this article. This study not only determines the analytical forms of optimal interval-weights, consistency index and consistency ratio, but also selects the best optimal weight set by introducing the secondary objective function. Some salient aspects of the proposed research are as follows. First, the mentioned analytic forms are obtained by developing a novel optimal modification based model, which is equivalent to the multiplicative model. This model provides a better understanding and interpretation of the underlying process of the multiplicative model. Second, the secondary objective function used for the selection of the best optimal weight set retains all the characteristics of the multiplicative model. Third, the analytic form of consistency ratio is useful in providing immediate feedback to the decision maker. Fourth, the analytic forms of optimal interval-weights and the best optimal weight set eliminate the need of an optimization software, reducing the computational complexity and enhancing the time efficiency.\\\\
	This research opens up an important future direction: can the proposed method be applied to obtain analytical forms of optimal weights for other models of BWM such as Euclidean BWM\cite{kocak2018euclidean}, Fuzzy BWM\cite{guo2017fuzzy}, etc., as well? And if not, how to derive an analytical framework for these models?
	\section*{Acknowledgements}
	The first author gratefully acknowledges the Council of Scientific \& Industrial Research (CSIR), India for financial support to carry out the research work. We would like to thank Dr. Jafar Rezaei for his useful comments that	improved the manuscript.

	\bibliographystyle{plain}
\end{document}